\documentclass[12pt, final]{amsproc}
\usepackage[cp850]{inputenc}
\usepackage{amssymb}
\usepackage[mathscr]{eucal}
\usepackage{amscd}
\usepackage{amsmath}
\usepackage{setspace}
\usepackage{tikz}
\def\Ker{\operatorname{Ker}}

\def\dim{\operatorname{dim}}

\theoremstyle{plain}
\newtheorem{theorem}{Theorem}
\newtheorem{proposition}{Proposition}
\newtheorem{lemma}{Lemma}
\newtheorem{corollary}{Corollary}
\newtheorem{definition}{Definition}
\theoremstyle{remark}
\newtheorem{remark}{Remark}

\title{A weak Hilbert space that is a twisted Hilbert space}

\address{Departamento de Matem\'aticas, Universidad de Extremadura, Avenida de Elvas s/n, 06011 Badajoz, Spain.}

\email{jesus@unex.es}
\author{Jes\'us Su\'arez de la Fuente}

\subjclass[2010]{(Primary) 46B20, 46B06; (Secondary) 46B70, 46M18, 46B45}

\keywords{Weak Hilbert, interpolation, twisted sum, centralizer}

\thanks{The author was supported in part by project MTM2016-76958-C2-1-P and project IB16056 of La Junta de Extremadura.}

\dedicatory{To Reyes and \'Ursula, my untwisted wings.}
\begin{document}

\begin{abstract}
We construct a weak Hilbert space that is a twisted Hilbert space. 
\end{abstract}

\maketitle
\section{Introduction}
The purpose of this paper is to construct, as announced, a weak Hilbert space that is a twisted Hilbert space. The new space $Z(T^2)$ is obtained as the derived space generated through complex interpolation between $T^2$, the 2-convexification  of the Tsirelson space, and its dual. Let $\mathcal Z$ denote either a Kalton-Peck space or an Enflo-Lindenstrauss-Pisier space (see below for definitions), then $\mathcal Z$ cannot be a subspace or a quotient of $Z(T^2)$ and $Z(T^2)$ is not isomorphic to either a subspace or a quotient of $\mathcal Z$. Roughly speaking, the only common subspaces or quotients of $Z(T^2$) and $\mathcal Z$ are the Hilbert spaces. 

Before stating in more detail the results herein, we recall some definitions and set our notation. In this paper we shall use the terminology commonly used in Banach space theory as it appears in \cite{AK}. ``Space" means ``infinite dimensional Banach space" unless specified otherwise. We will denote by $[x_j]_{j\in A}$ the closed linear span of a sequence $(x_j)_{j\in A}$ in a space $X$. If $X$ has an unconditional basis, and no confusion arises, $(e_j)_{j=1}^{\infty}$ will be such basis for which we also denote $X(A)=[e_j]_{j\in A}$ with $A \subseteq \mathbb N$. Given two spaces $X,Y$, we write $X\approx Y$ if they are isomorphic while for a finite $n$-dimensional Banach space $E$, we denote as usual
$$d_E=d(E, \ell_2^n),$$
where $d$ stands for the Banach-Mazur distance $d(X,Y)=\inf \{\|T\|\|T^{-1}\|\}$ and the infimum ranges over all isomorphims $T:X\to Y$.

Our construction is based on the Tsirelson space, for which the book of Casazza and Shura \cite{CS} provides a comprehensive study. To introduce Tsirelson's space we need to give inductively a sequence of norms. Pick an element $x$ in $c_{00}$, the vector space of all finitely supported sequences, and define:
$$\|x\|_0=\max |x_j|,$$
$$\|x\|_{m+1}= \max \left \{ \|x\|_m, \frac{1}{2} \max \left [ \sum_{j=1}^k \|A_j x\|   \right ] \right\} , \;\;(m\geq 0),$$
where the inner max is taken over all choices of finite subsets as $k$ varies and such that $k\leq A_1<...<A_k$. We have written $A_j<A_{j+1}$ as usual for $\max A_j <\min A_{j+1}$. Tsirelson's space $T$  is the completion of $c_{00}$ under the norm $\|x\|:=\lim_{m\to \infty}\|x\|_m$.  The $2$-convexification of $T$ is denoted $T^2$ and is the completion of $c_{00}$ under the norm
$$\left \|\sum x_j e_j \right \|_{T^2}=\left \| \sum |x_j|^2 e_j \right \|_T^{\frac{1}{2}}.$$
The main topic discussed here deals with the class of weak Hilbert spaces introduced by Pisier in \cite{Pi3}. 
\begin{definition}\label{weakhilbert}
We say $X$ is a weak Hilbert space if there is $0<\delta_0<1$ and a constant $C$ with the following property: every finite-dimensional subspace $E$ of $X$ contains a subspace $F\subseteq E$ with $\dim F\geq \delta_0 \dim E$ such that $d_F\leq C$ and there is a projection $P:X\to F$ with $\|P\|\leq C$. 
\end{definition}
The definition above is not the original one but is chosen out among the many equivalent characterizations given by Pisier \cite[Theorem 12.2.(iii)]{Pi}.  The canonical example of a weak Hilbert space different from the Hilbert space is $T^2$. Since the property ``to be a weak Hilbert space" passes to subspaces, quotients and duals (\cite{Pi}), we find that $(T^2)^*$ is also a weak Hilbert space.

To construct our example we make use of centralizers arising from interpolation so let us recall here some main points. Let $\omega$ denote the vector space of all complex scalar sequences with the topology of pointwise convergence and let $X$ be a space with an unconditional basis. An homogeneous map $\Omega:X\longrightarrow \omega$ is called a \textit{centralizer} if there is a constant $C$ so that: 
\begin{equation}\label{centralizer}
\| \Omega(ax)-a\Omega(x)\|_X\leq C\|a\|_{\infty}\|x\|_X,\;\;\:\;a\in \ell_{\infty},x\in X.
\end{equation}
A typical way to obtain centralizers is through complex interpolation so let $X_0,X_1$ be a couple of spaces with a joint unconditional basis and natural inclusions into $\omega$; the inclusions are not only linear but continuous by our choice of the topology on $\omega$. We shall consider the vector space $\mathcal F(X_0,X_1)$ of all functions $F:\mathbb S\to \omega$, which are bounded and continuous on the strip $$\mathbb S=\{ z: 0\leq Rez \leq 1 \},$$
and analytic on the open strip $\{z: 0< Re z <1 \},$ and moreover, the functions $t\in \mathbb R\to F(j+it)\in X_j$ with $j=0,1$ are bounded and continuous functions, which tend to zero as $|t|\to \infty$. $\mathcal F(X_0,X_1)$ is a Banach space when is endowed with the norm $$\|F\|_{\mathcal F}=\max \left ( \sup_{t\in \mathbb R} \| F(it)\|_{X_0},   \sup_{t\in \mathbb R} \| F(1+it)\|_{X_1}\right).$$
 The interpolation space $X_{\theta}:=(X_0,X_1)_{\theta}$ consists of all $x\in \omega$ such that $x=F(\theta)$ for some $F\in \mathcal F(X_0,X_1)$ endowed with the quotient norm $$\|x\|_{X_{\theta}}=\inf \{ \|F\|_{\mathcal F} : F(\theta)=x\}.$$
  We denote as usual by $\delta_{\theta}$ the natural quotient map $\mathcal F(X_0,X_1) \to  X_{\theta}$ with $0<\theta<1$.

Let us show how there is a natural centralizer on $X_{\theta}$ for which we assume for the sake of clarity $X_{\theta}$ has the same unconditional basis than $X_0,X_1$. First of all fix a constant $\rho>1$ and thus for each $x\in X_{\theta}$ pick a map $B(x)\in \mathcal F(X_0,X_1)$ with $B(x)(\theta)=x$ and $\|B(x)\|_{\mathcal F}\leq \rho\|x\|_{\theta}$. We shall refer to $B$ as a $\rho$-bounded selector for $\delta_{\theta}$ for which there is no loss of generality assuming that is homogeneous. A centralizer $\Omega_{\rho}:X_{\theta}\longrightarrow \omega$ comes defined as $$\Omega_{\rho}(x)=\delta'_{\theta}B(x).$$ 
 Centralizers and twisted sums are naturally related as we are about to see. Recall first that a \textit{short exact sequence} is a diagram like
\begin{equation*}\label{ses}
\begin{CD} 0 @>>>Y@>j>>Z@>q>>X@>>>0
\end{CD}
\end{equation*} where the morphisms are linear and continuous and such that the image of each arrow is the kernel of the next one. This condition implies that $Y$ is a subspace of $Z$ through $j$ and thanks to
the open mapping theorem we find that $X$ is isomorphic to $Z/j(Y)$. We refer to $Z$ as a \textit{twisted sum} of $Y$ and $X$ and if $Y\approx\ell_2\approx X$ we simply say that $Z$ is a \textit{twisted Hilbert space}. On the other hand, the centralizer $\Omega_{\rho}$ produces a natural twisted sum 
$$
\begin{CD} 0 @>>>X_{\theta}@>j>>d_{\Omega_{\rho}}X_{\theta}@>q>>X_{\theta}@>>>0,
\end{CD}
$$
 where
$$d_{\Omega_{\rho}}X_{\theta}:=\{  (x,y)\in \omega\times \omega: \|x-\Omega_{\rho}(y)\|+\|y\|<\infty \},$$ 
 endowed with the quasi-norm
\begin{equation}\label{quasinorm}
\|(x,y)\|=\|x-\Omega_{\rho}(y)\|+\|y\|.
\end{equation} 
The maps are simply given by $j(x)=(x,0)$ and $q(x,y)=y$. The quasi-norm (\ref{quasinorm}) is equivalent to the norm (\ref{norm}) of the space $dX_{\theta}$ called \textit{the derived space} and defined as the set of couples $(x,y)\in \omega\times \omega$ for which the following expression 
\begin{equation}\label{norm}
\|(x,y)\|_{dX_{\theta}}=\inf \{ \|F\|: F(\theta)=y, F'(\theta)=x\}
\end{equation}
is finite. The formal statement of this last is the following (see \cite[Proposition 7.2.]{CCS}).
\begin{proposition}\label{equivalence}
The space $dX_{\theta}$ and $d_{\Omega_{\rho}}X_{\theta}$ coincide, with equivalent quasi-norms.
\end{proposition}
Recall that a centralizer $\Omega_{\rho}$ is \textit{bounded} if $\Omega_{\rho}(y)\in X_{\theta}$  for every $y\in X_{\theta}$ and there is $K>0$ such that $\|\Omega_{\rho}(y)\|_{X_{\theta}}\leq K \|y\|_{X_{\theta}}$ in which case (\ref{quasinorm}) is equivalent to the product norm. It is clear that there is some ambiguity in the choice of $\Omega_{\rho}$ but the difference of any two choices is a bounded centralizer what implies that the quasi-norms (\ref{quasinorm}) are equivalent. One last important fact for us about $dX_{\theta}$ is the following.
\begin{proposition}\label{lambda} Let $(e_j)_{j=1}^{\infty}$ be the unconditional basis of $X_{\theta}$ and set $v_{2j-1}=(e_j,0), v_{2j}=(0,e_j)$ for $j\in \mathbb N$. Then
\begin{enumerate}
\item[$(i)$] The sequence $(v_j)_{j=1}^{\infty}$ is a basis for $dX_{\theta}$.
\item[$(ii)$] The sequence $(v_{2j})_{j=1}^{\infty}$ is unconditional.
\end{enumerate}
\end{proposition}
The first part is proved by adapting the proof of \cite[Theorem 4.10]{kaltpeck} while the second follows picking $a\in \{-1,1\}^{n}$ in (\ref{centralizer}).

Since we claim that our twisted-Hilbert is weak-Hilbert we should be able to distinguish it from the previous examples of twisted Hilbert spaces: The Kalton-Peck spaces and the Enflo-Lindenstrauss-Pisier example. Let us briefly outline these constructions. Given a Lipschitz map $\phi:[0,\infty)\to \mathbb R$ with $\phi(0)=0$, the Kalton-Peck space $Z_2(\phi)$ is the twisted Hilbert space $d_{\Omega}\ell_2$ generated with the centralizer 
\begin{eqnarray}\label{kaltpeck2}
\Omega(y)=\sum_{j=1}^{\infty} y_j \phi \left(\log \left( \|y\|_2/|y_j|\right)\right)e_j,
\end{eqnarray}
taking $0\cdot \phi(\log \infty)$ as $0$ and where $y=\sum y_je_j$. The case $\phi(t)=t$ corresponds to the so-called space $Z_2$. That $Z_2(\phi)$ contains the following sequence space with symmetric basis
\begin{equation}\label{orlicz} 
\ell_{M}:=[(0,e_j)]_{j=1}^{\infty}
\end{equation}
will be a crucial point to differentiate it from our example. The symmetry of the basis follows from Proposition \ref{lambda} and a quick look at (\ref{kaltpeck2}). We will also deal with the Enflo-Lindenstrauss-Pisier example given in \cite[Proposition 2]{ELP} where it is proved the existence for each $n\in\mathbb N$ of a twisted sum
$$
\begin{CD} 0 @>>>\ell_2^{n^2}@>>>\mathcal E_n@>>>\ell_2^n@>>>0,
\end{CD}
$$
for which the norm of any projection of $\mathcal E_n$ onto $\ell_2^{n^2}$ is lower bounded by $c(\log n)^{1/2}$, for some absolute constant $c>0$. Then the Enflo-Lindenstrauss-Pisier space is $\ell_2(\mathcal E_n)$.
\section{The example}\label{secbeta}
Let us introduce the space $Z(T^2)$. It is clear that $\ell_2$ is continuously and densely embedded in $T^2$ and therefore $(T^2,(T^2)^*)_{\frac{1}{2}}=\ell_2$ with equivalent norms by \cite[Corollary 4.3]{CoS}. We define $Z(T^2):=d(T^2,(T^2)^*)_{\frac{1}{2}}$ and the induced short exact sequence
\begin{equation*}
\begin{CD} 0 @>>>\ell_2@>>>Z(T^2)@>>>\ell_2@>>>0
\end{CD}
\end{equation*}
shows that $Z(T^2)$ is a twisted Hilbert space. That $Z(T^2)$ is not Hilbert can be seen as follows: it does not have type $2$ since otherwise the unconditionality of $(v_{2j})_{j=1}^{\infty}$ would imply that the centralizer providing an equivalent quasi-norm is bounded which, by  Kalton's theorem \cite[Theorem 7.6.]{ka}, implies $T_2\approx \ell_2$, which is absurd.

The key step to show that  $Z(T^2)$ is a weak Hilbert space is to compute the Rademacher type $2$ constants $a_{n,2}(X)$ for certain $\log(\log(n))$-codimensional subspaces of $Z(T^2)$. Recall that for a space $X$ the number $a_{n,2}(X)$ is defined to be the least constant $a$ such that
$$\mathbb E\left \|\sum_{j=1}^n \varepsilon_j x_j\right\|\leq a \left ( \sum_{j=1}^n \| x_j\|^2\right)^{1/2},$$
for all $x_1,...,x_n\in X$ and where the average is taken over all $(\varepsilon_j)_{j=1}^n\in \{-1,1\}^n$. The space $X$ has \textit{type $2$} if $a_2(X):=\sup_{n\in \mathbb N}a_{n,2}(X)<\infty$. The cotype 2 constant $c_{n,2}(X)$ is defined in a similar vein as the least constant $c$ such that $$\left ( \sum_{j=1}^n \| x_j\|^2\right)^{1/2}\leq c\cdot \mathbb E\left \|\sum_{j=1}^n \varepsilon_j x_j\right\|,$$
for all $x_1,...,x_n\in X$ and thus $X$ has \textit{cotype $2$} if $c_2(X):=\sup_{n\in \mathbb N}c_{n,2}(X)<\infty$. It will be very useful for us to introduce another norm on $\mathcal F(X_0,X_1)$ given in \cite{Cald}. For a fixed $0<\theta<1$ and given $F\in\mathcal F(X_0,X_1)$, we define
$$\|F\|_{\mathcal F_{\theta}}:=\int_{-\infty}^{\infty} \|F(it)\|_{X_0}dP_0(\theta,t)+\int_{-\infty}^{\infty} \|F(1+it)\|_{X_1}dP_1(\theta,t),$$
where $P_0,P_1$ denote the Poisson kernerls on $\mathbb S$ (see \cite[Page 93]{BL} for the definition). Let us observe that for each $F\in \mathcal F(X_0,X_1)$,
\begin{equation}\label{calderon}
\|F(\theta)\|_{X_{\theta}}\leq \|F\|_{\mathcal F_{\theta}}\leq \|F\|_{\mathcal F}.
\end{equation}
The first inequality is a simple application of \cite[Lemma 4.3.2.(iii)]{BL} while the second is due to the well known fact that $\int_{-\infty}^{\infty}dP_0(\theta,t)=1-\theta$  and $\int_{-\infty}^{\infty}dP_1(\theta,t)=\theta$. We follow now an idea of a recent paper of Castillo, Ferenczi and Gonz\'alez \cite{CFG} where they introduced an useful inequality. We prove an ``average version" in the spirit of \cite[Lemma 4.8.]{CFG}.
\begin{lemma}\label{yo}
Let $X_0,X_1$ be spaces with a joint unconditional basis and let $\Omega_{\rho}$ be a centralizer induced by $X_{\theta}=(X_0,X_1)_{\theta}$ with $0<\theta<1$. Given $b_1,...,b_n\in X_{\theta}$ we have
\begin{equation}\label{yo2}
\mathbb E\left \| \Omega_{\rho} \left(\sum_{j=1}^n \varepsilon_j b_j \right)- \sum_{j=1}^n\varepsilon_j \Omega_{\rho} (b_j)- \log \frac{a_{n,2}(X_0)}{a_{n,2}(X_1)}  \sum_{j=1}^n \varepsilon_j b_j \right\|_{X_{\theta}}\leq \gamma \cdot \left(\sum_{j=1}^n \|b_j\|_{X_{\theta}}^2\right)^{1/2},
\end{equation}
where $$\gamma= \widetilde{\theta}\cdot \rho \cdot a_{n,2}(X_0)^{1-\theta}a_{n,2}(X_1)^{\theta}.$$
\end{lemma}
\begin{proof} Let $B$ be a $\rho$-bounded selector for the quotient map $\delta_{\theta}$ so that $\Omega_{\rho}=\delta_{\theta}'B$. Fix $b_j\in X_{\theta}$ for $j\leq n$ and given $\varepsilon=(\varepsilon_j)_{j=1}^n\in \{-1,1\}^n$ define $$F_{\varepsilon}(z)=\frac{\sum_{j=1}^n \varepsilon_j B(b_j)(z)}{a_{n,2}(X_0)^{1-z}a_{n,2}(X_1)^z}\;\;\text{with}\;z\in \mathbb S.$$
Since every $F_{\varepsilon}$ is, up to a bounded scalar function, a linear combination of elements $B(b_j)\in \mathcal F(X_0,X_1)$, it is a straightforward computation that $F_{\varepsilon}\in \mathcal F(X_0,X_1)$. Observe that 
$$a_{n,2}(X_0)^{1-\theta}a_{n,2}(X_1)^{\theta}\delta'_{\theta}(F_{\varepsilon}(z))=\sum_{j=1}^n\varepsilon_j\Omega_{\rho}(b_j)+\log \frac{a_{n,2}(X_0)}{a_{n,2}(X_1)}\sum_{j=1}^n \varepsilon_j b_j.$$
Therefore, we have 
\begin{eqnarray*}
\Lambda_{\varepsilon}&:=& \left\| \Omega_{\rho} \left(\sum_{j=1}^n \varepsilon_j b_j \right)- \sum_{j=1}^n\varepsilon_j \Omega_{\rho} (b_j)- \log \frac{a_{n,2}(X_0)}{a_{n,2}(X_1)}  \sum_{j=1}^n \varepsilon_j b_j \right\| \\
&=&\left\|\Omega_{\rho}\left(\sum_{j=1}^n \varepsilon_j b_j\right)-\delta'_{\theta}(a_{n,2}(X_0)^{1-\theta}a_{n,2}(X_1)^{\theta}F_{\varepsilon})\right \|\\
&=&\left\| \delta'_{\theta}\left (B\left( \sum_{j=1}^n\varepsilon_j b_j\right)-a_{n,2}(X_0)^{1-\theta}a_{n,2}(X_1)^{\theta}F_{\varepsilon} \right) \right\|.
\end{eqnarray*}
It is clear that 
\begin{equation}\label{kernel}
B \left( \sum_{j=1}^n \varepsilon_j b_j\right)-a_{n,2}(X_0)^{1-\theta}a_{n,2}(X_1)^{\theta}F_{\varepsilon}\in \Ker \delta_{\theta}.
\end{equation}
A quick check will show that given $F\in \Ker \delta_{\theta}$, the map $G(z):=(z-\theta)^{-1}F(z)\in \mathcal F(X_0,X_1)$ satisfies $F'(\theta)=G(\theta)$ (\cite{CFG}), so we find by (\ref{calderon}) that $\|F'(\theta)\|_{X_{\theta}}\leq \|G\|_{\mathcal F_{\theta}}$. And then a trivial computation gives $ \|G\|_{\mathcal F_{\theta}}\leq \widetilde{\theta}\cdot \|F\|_{\mathcal F_{\theta}}$ for $\widetilde{\theta}=\max\{ \theta^{-1}, (1-\theta)^{-1} \}$. In particular, if we apply this argument to the map in (\ref{kernel}) we have
\begin{eqnarray*}
\Lambda_{\varepsilon}&\leq&  \widetilde{\theta}  \left (\left\|B\left (\sum_{j=1}^n\varepsilon_jb_j\right)\right\|_{\mathcal F_{\theta}} +\left \| a_{n,2}(X_0)^{1-\theta}a_{n,2}(X_1)^{\theta}F_{\varepsilon}\right\|_{\mathcal F_{\theta}}\right)\\
\end{eqnarray*} and averaging, using also the second inequality of (\ref{calderon}) for the first term, we find
\begin{eqnarray}\label{ave}
\mathbb E \Lambda_{\varepsilon} \leq \widetilde{\theta}\cdot \rho \cdot \mathbb E \left \| \sum_{j=1}^n\varepsilon_jb_j\right\| + \widetilde{\theta}\cdot a_{n,2}(X_0)^{1-\theta}a_{n,2}(X_1)^{\theta}\cdot\mathbb E \|F_{\varepsilon}\|_{\mathcal F_{\theta}} .
\end{eqnarray}
We have the following bound for the second term on the right side of (\ref{ave})
\begin{eqnarray*}
\mathbb E \| F_{\varepsilon}\|_{\mathcal F_{\theta}}&=& \sum_{\varepsilon\in \{-1,1\}^n} \frac{1}{2^n} \left\|  \frac{\sum_{j=1}^n \varepsilon_j B(b_j)(z)}{a_{n,2}(X_0)^{1-z}a_{n,2}(X_1)^z}  \right\|_{\mathcal F_{\theta}}\\
&=& \sum_{k=0}^1\int_{-\infty}^{\infty} \sum_{\varepsilon\in \{-1,1\}^n} \frac{1}{2^n} \left \|\frac{\sum_{j=1}^n \varepsilon_j B(b_j)(k+it)}{a_{n,2}(X_0)^{1-k-it}a_{n,2}(X_1)^{k+it}}\right  \|_{X_k} dP_k(\theta,t)\\
&\leq&\sum_{k=0}^1\int_{-\infty}^{\infty} \left( \sum_{j=1}^n  \left \| B(b_j)(k+it)\right  \|_{X_k}^2 \right)^{1/2} dP_k(\theta,t)\\
&\leq& \sum_{k=0}^1\int_{-\infty}^{\infty} \left( \sum_{j=1}^n  \left \| B(b_j)\right  \|_{\mathcal F}^2 \right)^{1/2} dP_k(\theta,t)\\
&\leq & \rho\cdot \sum_{k=0}^1\int_{-\infty}^{\infty} \left( \sum_{j=1}^n  \left \| b_j\right  \|_{X_\theta}^2 \right)^{1/2} dP_k(\theta,t)=\rho\cdot \left( \sum_{j=1}^n  \left \| b_j\right  \|_{X_\theta}^2 \right)^{1/2}.
\end{eqnarray*}
We may trivially bound the first term on the right side of (\ref{ave}) as $$\mathbb E \left \| \sum_{j=1}^n\varepsilon_jb_j\right\|\leq a_{n,2}(X_{\theta})\left( \sum_{j=1}^n \|b_j\|^2\right)^{1/2},$$ so it only remains to check that $a_{n,2}(X_{\theta})\leq a_{n,2}(X_0)^{1-\theta}a_{n,2}(X_1)^{\theta}$ to finish the proof. This follows by a typical log-convexity argument. Indeed, pick $b_1,...,b_n\in X_{\theta}$ and given $\delta>0$ find $F_j\in \mathcal F$ such that $F_j(\theta)=b_j$ and $\|F_j\|\leq (1+\delta)\|b_j\|$ for $j\leq n$. Apply the bound given in \cite[Lemma 4.3.2.(ii)]{BL} to $\left\|\sum_{j=1}^n\varepsilon_j F_j(\theta)\right\|$ and average, then use H\"older's inequality and our choice of $F_j$ to conclude letting $\delta\to 0$.
\end{proof}
The idea of averaging in the inequality \cite[Lemma 4.8.]{CFG} was introduced in \cite{Co} under some restrictions for analytic families of Banach spaces.  Lemma \ref{yo} admits an analogous formulation in terms of the constants $a_{n,p}(X)$ and the proof also recovers the case of \cite[Lemma 2.11.]{Co} although we set no restrictions on the type of the spaces involved. Getting back to our business, we are interested in computing $a_{n,2}(Z(T^2))$. 
\begin{proposition}\label{tipoder} Let $X_0,X_1$ be spaces with a joint unconditional basis and let $\Omega_{\rho}$ be a centralizer induced by $X_{\theta}=(X_0,X_1)_{\theta}$ with $0<\theta<1$. Then
$$a_{n,2}(d_{\Omega_{\rho}}X_{\theta})\leq \widetilde{\theta}\cdot\rho\cdot a_{n,2}(X_{\theta})\left(a_{n,2}(X_0)^{1-\theta}a_{n,2}(X_1)^{\theta} + \left|\log \frac{a_{n,2}(X_0)}{a_{n,2}(X_1)}\right|\right).$$
\end{proposition}
\begin{proof}
Pick vectors $x_j=(a_j,b_j)$ in $d_{\Omega_{\rho}}X_{\theta}$ for $j\leq n$ and let us denote by simplicity
$$\Delta_{\varepsilon}=\left \|\Omega_{\rho} \left(\sum_{j=1}^n \varepsilon_j b_j\right)- \sum_{j=1}^n \varepsilon_j\Omega_{\rho}(b_j)\right\|,\;\; \varepsilon\in \{-1,1\}^n.$$ Then 
\begin{eqnarray*}
\mathbb E\left \| \sum_{j=1}^n \varepsilon_j x_j \right \| &\leq & \mathbb E \left  \| \sum_{j=1}^n \varepsilon_j(a_j - \Omega_{\rho} (b_j))\right\| + \mathbb E\left \|  \sum_{j=1}^n \varepsilon_j b_j \right \|+ \mathbb E  \Delta_{\varepsilon}\\
&\leq & 2a_{n,2}(X_{\theta})\left(  \sum_{j=1}^n \|x_j\|^2\right) ^{1/2}+\mathbb E  \Delta_{\varepsilon} .\\
\end{eqnarray*}
A quick computation finishes the proof since we have  by Lemma \ref{yo} $$\mathbb E  \Delta_{\varepsilon} \leq \gamma \cdot\left(\sum_{j=1}^n\|b_j\|^2\right)^{1/2}+\mathbb E\left\|\log \frac{a_{n,2}(X_0)}{a_{n,2}(X_1)} \sum_{j=1}^n b_j\right\|.$$
\end{proof}
We are about to state and prove that $Z(T^2)=[v_j]_{j=1}^{\infty}$ is weak-Hilbert so we are interested in the $2(n-1)$-codimensional subspaces $V_{n}:=[v_j]_{j=2n-1}^{\infty}$. Let us fix some notation for the rest of the paper to simplify the forthcoming arguments. Let $(e_j)_{j=1}^{\infty}$ and $(e_j^*)_{j=1}^{\infty}$ be the basis in $T^2$ and $(T^2)^*$ respectively and let us write $E_n=[e_j]_{j=n}^{\infty}$ and $E_n^*=[e_j^*]_{j=n}^{\infty}$. Once this is done, it is not hard to see that $$V_{n}=d(E_n,E_n^*)_{\frac{1}{2}}.$$ 
It is important to recall that $Z(T^2)$ is $\lambda$-isomorphic to its own dual where the duality comes defined, given $(x,y)\in Z(T^2)^*$ and $(u,v)\in Z(T^2)$, by the rule
\begin{equation}\label{duality}
\langle (x,y),(u,v)  \rangle= \langle x,v\rangle + \langle y,u\rangle,
\end{equation}
 see \cite[Corollary 4]{FC}. It is therefore plain to deduce that $V_n$ is also $\lambda$-isomorphic to its own dual for every $n\in \mathbb N$.
The following hierarchy of functions will be useful from now on. 
\[  \left\{ \begin{array}{ll}
      \displaystyle \log_1(n)=\log(n), & \mbox{and}\\
        \log_{m+1}(n)=\log(\log_m(n)) & \mbox{for $n$ large enough so that $\log_{m+1}(n)\geq 1$, ($m\geq 1$).}
        \end{array} \right. \] 
We state and prove the main result of the paper.
\begin{theorem}\label{oneyear}
$Z(T^2)$ is a weak Hilbert space.
\end{theorem}
\begin{proof} Our main tool to prove the claim is a lemma of Johnson \cite[Lemma 1.6]{Jo} showing that if the $5^{(5^n)}$-dimensional subspaces of certain $n$-codimensional subspaces of $X$ are euclidean then $X$ is weak-Hilbert. Let us divide the proof in several steps for the sake of clarity.
\begin{flushleft}\textbf{Step I.} \textit{ There is $C>0$ such that if $E$ is a $5^{(5^n)}$-dimensional subspace of $V_{n}$ then  $d_E\leq C$.}
\end{flushleft}
 Indeed, we will prove a stronger statement, namely that given $m=1,2,...$, there is some $C_m>0$ such that every $n$-dimensional subspace $E$ of $V_{\log_m(n)}$ satisfies $d_E\leq C_m$, what in particular implies Step I for $m\geq 2$. Thus, fix $m=1,2,...$, and pick $k=\log_m(n)$. We show now that there is $c_m>0$ such that for every $n\in\mathbb N$ $$\max\left(a_{n,2}(E_k),a_{n,2}(E_k^*) \right)\leq c_m.$$ Indeed, $T^2$ has type 2 so it follows $\sup_{n\in \mathbb N}a_{n,2}(E_k)\leq a_2(T^2)<\infty$. \cite[Proposition Ab2]{CS} also holds in $(T^2)^*$ so there is $c>0$ such that for every subspace $E$ in $E_k^*$ with $\dim E \leq n$ we have $d_E\leq c^m$. Since we always have
\begin{equation}\label{typedist} 
\mathbb E\left \|  \sum_{j=1}^n \varepsilon_j x_j\right\|\leq d_{[x_j]_{j=1}^n}\left( \sum_{j=1}^n \|x_j\|^2 \right)^{1/2},
\end{equation}
we find that $a_{n,2}(E_k^*)\leq c^m$ so that $c_m=\max\{a_2(T^2),c^m\}$ works. Now, since $V_{k}=d(E_k,E_k^*)_{\frac{1}{2}}$, we find by Proposition \ref{tipoder} and Proposition \ref{equivalence} that
$$a_{n,2}(V_{k})\leq f(c_m),$$ for some function $f$. We are interested also in the cotype $2$ constants for which we have that $c_{n,2}(X^*)\leq a_{n,2}(X)$ for any space $X$, see the proof of \cite[Proposition 3.2.]{Pi2}. Since $V_{k}$ is $\lambda$-isomorphic to its dual, we find that
$c_{n,2}(V_{k})\leq \lambda\cdot a_{n,2}(V_{k})$ for some $\lambda>0$ independent of $n$. Pick $E$ a subspace of $V_{k}$ with $\dim E =n$, then
\begin{eqnarray*}
d_E&\leq& a_2(E) \cdot  c_2(E)\\
&\leq & 2\sqrt{2\pi}  \cdot a_{n,2}(E)\cdot c_{n,2}(E)\\
&\leq &  2\sqrt{2\pi} \cdot  a_{n,2}(V_{k})\cdot  c_{n,2}(V_{k})\\
&\leq & 2\sqrt{2\pi} \cdot \lambda \cdot a_{n,2}(V_{k})^2\\
&= & 2\sqrt{2\pi} \cdot \lambda \cdot f(c_m)^2,
\end{eqnarray*}
where the first inequality is a well known result of Kwapie\'n \cite{Kw} and the second holds by a result of Tomczak-Jaegermann \cite[Theorem 2]{To} used in the form $a_2(E)\leq \sqrt {2\pi}\cdot a_{n,2}(E)$ and $c_2(E)\leq 2\cdot c_{n,2}(E)$. Therefore, set $C_m=2\sqrt{2\pi} \cdot \lambda \cdot f(c_m)^2$.\\

\begin{flushleft}\textbf{Step II.}  \textit{There is $C>0$ such that if $E$ is an $n$-dimensional subspace of $V_n$ then $d_E\leq C$ and there is a projection $P:Z(T^2)\to E$ with $\|P\|\leq C$.}
\end{flushleft}
By Step I, it follows from the already mentioned result of Johnson \cite[Lemma 1.6]{Jo} that every $5^n$-dimensional space of $V_{n}^*$ is $3C$-isomorphic to a Hilbert space and $3C$-complemented in $V_{n}^*$. Therefore, again by \cite[Lemma 1.6]{Jo},  every $n$-dimensional subspace of $V_{n}$ is $9C$-isomorphic to a Hilbert space and $9C$-complemented in $V_{n}$, thus $9C(K+1)$-complemented in $Z(T^2)$ where $K$ denotes the basis constant of $(v_j)_{j=1}^{\infty}$.\\

 \begin{flushleft} \textbf{Step III.} \textit{$Z(T^2)$ is a weak Hilbert space.}
\end{flushleft}
If $E$ is a $4n$-dimensional subspace of $Z(T^2)$, let $F:=E \cap [v_j]_{j=2n-1}^{\infty}$. Then $\dim F \geq 2n-1\geq n$ and $F$ is a subspace of  $V_n=[v_j]_{j=2n-1}^{\infty}$. Therefore $F$ (and thus $E$) contains an $n$-dimensional subspace that must be $C$-isomorphic to a Hilbert space and $C$-complemented in $Z(T^2)$ by Step II. If $E$ is an arbitrary $N$-dimensional subspace of $Z(T^2)$, with $N=4n+k$ and $0\leq k\leq 3$, it follows from the previous argument that $E$ contains $F$ with $\dim F=n\geq \frac{1}{7}(4n+k)=\frac{1}{7}\dim E$ that is $C$-complemented in $Z(T^2)$ and $d_F\leq C$. Thus, $Z(T^2)$ satisfies Definition \ref{weakhilbert} and therefore is a weak Hilbert space.
\end{proof}
We are in position now to distinguish $Z(T^2)$ from previous examples of twisted Hilbert spaces.
\begin{proposition}\label{notisom}
Let $\phi$ be a Lipschitz map such that either $\lim_{t \to \infty} \phi'(t)=0$ monotonically or $\phi(t)=ct$ for $c\neq 0$. Then
\begin{enumerate}
\item[$(i)$] $Z(T^2)$ is not isomorphic to either a subspace or a quotient of $Z_2(\phi)$. 
\item[$(ii)$] $Z(T^2)$ is not isomorphic to either a subspace or a quotient of $\ell_2(F_n)$ with $\dim F_n<\infty$ for all $n\in \mathbb N$. 
\end{enumerate}
\end{proposition}
\begin{proof}
Let us prove $(i)$. Suppose that $S:Z(T^2) \to  Z_2(\phi)$ is an into isomorphism and write $y_j=S(v_{2j})/\|S(v_{2j})\|$ for $j\in\mathbb N$. By \cite[Theorem 5.4.]{kaltpeck}, there is an infinite subset $A\subset \mathbb N$ such that the basic sequence $(y_j)_{j\in A}$ is equivalent to either $\ell_2$ or the usual basis of $\ell_{M}$ (see (\ref{orlicz})). Thus, we have two possibilities. Assume $(y_j)_{j\in A}$ is equivalent to the basis of $\ell_2$. Then $(v_{2j})_{ j\in A}$ in $Z(T^2)$ is also equivalent to $\ell_2$ so that if $\Omega$ denotes a centralizer generating $Z(T^2)$ we find that for some $K>0$ and every $(\lambda_j)_{j\in A}$ $$\left\| \sum_{j\in A }\lambda_j v_{2j}\right\|=\left \|\Omega \left(\sum_{j\in A }\lambda_j e_j\right)\right \|+ \left \| \sum_{j\in A}\lambda_j e_j\right \|\leq K \left\| \sum_{j\in A}\lambda_j e_j\right \|.$$ In other words, the centralizer $\Omega_{|A}$ is bounded. Let us recall that $(T^2(A), (T^2)^*(A))_{1/2}=\ell_2(A)$ (see \cite[Corollary 4.3]{CoS}) and we may trivially assume that $\Omega_{|A}$ is the corresponding centralizer. Since $\Omega_{|A}$ is bounded, we find that $T^2(A)\approx\ell_2(A)$  by \cite[Theorem 7.6.]{ka} which is impossible since $T^2$ contains no copies of $\ell_2$ \cite{CS}. Let us assume now that $(y_j)_{j\in A}$ is equivalent to $\ell_{M}$. Then, as before, $(v_{2j})_{ j\in A}$ in $Z(T^2)$ is equivalent to the symmetric basis of $\ell_{M}$ that must be a weak Hilbert space by Theorem \ref{oneyear}. Therefore, by \cite[Proposition 12.4.]{Pi}, $\ell_{M}$ and $\ell_2$ have equivalent norms and we are in the previous case.

For the second part of $(i)$ let us recall that given a centralizer on $\ell_2$, the corresponding twisted-Hilbert is isomorphic to its own dual (\cite[Corollary 4]{FC}). Therefore, a simple duality argument combined with the previous proof shows the claim.

To prove $(ii)$ let us suppose that $Z(T^2)$ is a subspace of $\ell_2(F_n)$ with $\dim F_n<\infty$. Recall that every subspace of $\ell_2(F_n)$ with $\dim F_n<\infty$ has the property that every weak null normalized sequence admits a subsequence equivalent to the unit vector basis of $\ell_2$ (see \cite[Corollary 4.6.]{OS} for a stronger statement). The normalized sequence $(v_{2j})_{j=1}^{\infty}$ in $Z(T^2)$ admits no subsequence equivalent to the unit vector basis in $\ell_2$ as we proved above. It is weakly null by (\ref{duality}) so we reach a contradiction. If $Z(T^2)$ is a quotient of $\ell_2(F_n)$ then, again by \cite[Corollary 4]{FC}, it follows that $Z(T^2)$ is a subspace of  $\ell_2(F_n^*)$ and we are in the previous case.
\end{proof}

It follows from the arguments contained in the previous proof that $Z_2(\phi)$ is not weak-Hilbert unless it is a Hilbert space (see also \cite[Proposition 4.1.]{M}). Hence, $Z_2$ is not isomorphic to a subspace or a quotient of $Z(T^2)$. All this suggests that the only common subspaces of $Z_2$ and $Z(T^2)$ are the Hilbert spaces. The proof of Theorem \ref{oneyear} has also the following consequence.
\begin{corollary}\label{oneyearplus}
There is $C>0$ such that for all finite-dimensional subspaces $E$ of $Z(T^2)$ and for all $m=1,2,...$, we have 
$$d_E\leq C^m\log_m(\dim E).$$
\end{corollary}
\begin{proof} Given a finite $n$-dimensional subspace $E$ of $V_1=Z(T^2)$, we find as in the last part of Step I in Theorem \ref{oneyear} the bound $d_E\leq 2\sqrt{2\pi} \cdot \lambda \cdot a_{n,2}(V_{1})^2$. Then the claim follows from this once we prove
\begin{eqnarray}\label{tipo}
\left[a_{n,2}(Z(T^2))\right]^2\leq C^m\log_m(n),
\end{eqnarray}
for some absolute constant $C$. First, let us compute the type 2 constants of $T^2$ and $(T^2)^*$. Since $T^2$ has type 2 (\cite[Theorem 13.1.]{Pi}), we find that $\sup_{n\in \mathbb N}a_{n,2}(T^2)=a_2(T^2)<\infty$. For $(T^2)^*$, recall that by \cite[Theorem Ab7]{CS} we have  $d_E=o(\log_m(\dim E))$ for every subspace $E$ of $(T^2)^*$. Indeed, a combination of the arguments in \cite[Theorem Ab5, Theorem Ab6]{CS} gives that $d_E\leq c^m\log_m(\dim E)$ where $c$ is a universal constant. And thus by (\ref{typedist}), since $\dim [x_j]_{j=1}^n\leq n$, we have that $a_{n,2}((T^2)^*)\leq c^m\log_m(n)$. Therefore, by Proposition \ref{tipoder} and Proposition \ref{equivalence}, we trivially find (\ref{tipo}).
\end{proof}
In particular, if a sequence $(F_n)_{n=1}^{\infty}$ of finite-dimensional spaces are uniformly isomorphic to subspaces of $Z(T^2)$ then for every $m=1,2,...,$ it must be
\begin{equation}\label{estimate}
d_{F_n}=o(\log_m(\dim F_n)).
\end{equation}
If we denote $Z_2^n=[v_j]_{j=1}^{2n}$ endowed with the norm of $Z_2$ then from \cite[Theorem 6.3.]{kaltpeck} we infer that there is $c>0$ such that $d_{Z_2^n}\geq c\log(n)$ so that $Z_2^n$ does not satisfy (\ref{estimate}). Since an analogous result holds for $\mathcal E_n$ by the results in \cite{ELP}, we also find the following.
\begin{corollary}
$\ell_2(\mathcal E_n)$ is not isomorphic to a subspace or a quotient of $Z(T^2)$.
\end{corollary}
Let us recall that the first example of weak-Hilbert that is $\ell_2$-saturated was given by Edgington \cite{Ed}, while the first example with no unconditional basis was introduced by Komorowski \cite{Kom}. Later more examples were given by Komorowski and Tomczak-Jaegermann \cite{KomT}; see also \cite{ACK}. We have
\begin{proposition} 
$Z(T^2)$ has no unconditional basis and is $\ell_2$-saturated.
\end{proposition}
\begin{proof}
To prove the first claim let us assume $Z(T^2)$ has an unconditional basis. Then by \cite[Theorem 2.3.]{ka2} we have $Z(T^2)\approx \ell_2$ that is absurd. The second claim holds because every twisted Hilbert space is $\ell_2$-saturated \cite{CG}.
\end{proof}
Although $Z(T^2)$ has no unconditional basis, it is not hard to see that $Z(T^2)$ has an unconditional decomposition into two-dimensional subspaces. In general, a result of Maurey and Pisier included in \cite{M2} shows that weak Hilbert spaces enjoy FDD. 

Recall that a space $X$ has the \textit{Maurey extension property} if for every subspace $X_0$ of $X$, every operator from $X_0$ in a cotype 2 space admits an extension to the whole $X$ (\cite[Definition 2.1.(i)]{CaNi}). It follows from results of Milman and Pisier \cite[Theorem 10, Remark 11]{MiPi} that a space with the Maurey extension property has weak type 2 while the converse is not known, see \cite[pages 154, 155]{CS} and \cite[p. 2127]{CP}. We  observe that such a converse does not hold.
\begin{proposition}
$Z(T^2)$ has weak type 2 and fails the Maurey extension property.
\end{proposition}
\begin{proof}
$Z(T^2)$ has weak type 2 by the very definition of a weak Hilbert space (see \cite[Definition 12.1.]{Pi}) and fails the Maurey extension property since it contains an uncomplemented copy of $\ell_2$.
\end{proof}
 Indeed, we have proved that $Z(T^2)$ has not the \textit{Maurey projection property} (see \cite[Definition Aa11]{CS}) nor the property $M_2$ (see \cite[Definition 2.1.(ii)]{CaNi} and also \cite[Theorem 2.3.]{CaNi}). 
\begin{remark}
It is known that weak Hilbert spaces are \textit{ergodic} \cite[Corollary 2.4.]{An}. Therefore, $Z(T^2)$ is an example of an ergodic twisted Hilbert space not isomorphic to $\ell_2$. This answers a question raised in \cite{Cu}.
\end{remark}

{\bf Acknowledgements.}\label{ackref}
Part of this work in an early stage was done while the author was visiting William B. Johnson at Texas A\&M University around the SUMIRFAS conference 2017. We would like to thank him for the hospitality. We are grateful to Jes\'us M. F. Castillo for his interest in this work. We are also deeply indebted with the referee for the careful reading of the manuscript and the insightful suggestions.

\end{document}